
\magnification=\magstep1
\input amstex
\documentstyle{amsppt}
\hsize = 6.5 truein
\vsize = 9 truein
\TagsAsMath
\NoRunningHeads
\topmatter
\title
On G\"odel's Second Incompleteness Theorem
\endtitle
\author
Thomas Jech
\endauthor
\thanks Supported in part by NSF grant DMS--8918299
\endthanks
\address 
Department of Mathematics,
The Pennsylvania State University,
University Park, PA 16802
\endaddress
\email jech\@math.psu.edu 
\endemail
\endtopmatter
 
\document
\baselineskip 18pt
 
G\"odel's Second Incompleteness Theorem states that no
sufficiently strong consistent mathematical theory can prove its own
consistency. In \cite{1} this is proved for every axiomatic
theory extending the Peano Arithmetic. For axiomatic set theory
a simpler proof was given in \cite{2} using the fact that in
set theory, consistency of a set of axioms is equivalent to the
existence of a model.

In this note we give a very simple proof of G\"odel's Theorem for
set theory: 

\proclaim {Theorem} 
It is unprovable in set theory (unless it is inconsistent) that
there exists a model of set theory.
\endproclaim

By ``set theory'' we mean any axiomatic set theory  with finitely
many axioms. The
proof (suitably modified) works for any theory
sufficiently strong to formulate the
concepts ``model'', ``satisfies'' and ``isomorphic'',
such as second order arithmetic and its weaker versions.

If $M$ and $N$ are models of set theory (henceforth {\it models}),
we define
$$M<N \text{ if there exists some } M' \in N \text{ isomorhic to } M.$$
If $M<N$ and if $M'\in N$ is isomorphic to $M$ then for every sentence
$\sigma,$
$$ 
M \models \sigma \quad\text{ if and only if } \quad N \models 
(M'\models\sigma).\tag1
$$
It follows that
$$
\text{ if } M_1<M_2 \text{ and } M_2<M_3 \text{ then } M_1<M_3. \tag2
$$

Let us consider some fixed coding of formulas by numbers (G\"odel
numbering), and let $S_n$ be the name for the $n$th definable
set of numbers.

{\bf Definition. } $S$ is the set of all numbers $n$ with the property
that for every model $M$ there is a model $N<M$ such that 
$N \models n \notin S_n.$
\medpagebreak

Let $k$ be the G\"odel number of $S.$ The following sentence is
provable in set theory:
$$
k \in S_k \,\leftrightarrow \,\forall M \,\exists N<M \, 
N\models (k \notin S_k).
$$
Therefore (by (1) and (2)), if $M$ is any model then
$$
M \models (k\in S_k) \,\leftrightarrow \,\forall M_1 < M \,\exists M_2 < M_1 \,
M_2 \models (k\notin S_k).\tag3
$$

We shall now assume that it is provable that there is a model, and prove
a contradiction. As a consequence of the assumption, we have not only
$$
\text{there exists a model}\tag4
$$
but also (using (1))
$$
\text{for every model } M \text{ there exists a model } N < M. \tag5
$$

Toward a contradiction, we assume first that $k \in S.$ Let $M_1$
be an arbitrary model (by (4)); there exists an $M_2 < M_1$ that 
satisfies $k \notin S_k.$ By (3) there exists an $M_3 < M_2$
such that for every $M_4 < M_3,$ $M_4 \models k \in S_k.$ Therefore
$k\notin S,$ a contradiction.

Thus assume that $k \notin S.$ There is an $M_1$ such that for every
$M_2 < M_1,$ $M_2 \models k \in S_k.$ Let $M_2$ be an arbitrary model
$< M_1$ (it exists by (5)). As $M_2 \models k \in S_k,$ we have (by (3))
$$
\forall M_3 < M_2 \,\exists M_4 < M_3 \, M_4 \models(k \notin S_k).
$$
Let $M_3 < M_2$ (by (5)) and let $M_4 < M_3$ be such that
$M_4 \models k \notin S_k.$ As $M_4 < M_1$ (by (2)), we have
a contradiction.

\enddocument
 
\bye